\documentclass[11pt]{amsart}
\usepackage{amssymb, latexsym}
\theoremstyle{plain}
\newtheorem{theorem}{Theorem}

\newtheorem{proposition}{Proposition}

\theoremstyle{remark}

\newtheorem*{Remark 1}{Remark 1}
\newtheorem*{Remark 2}{Remark 2}
\newtheorem*{Remark 3}{Remark 3}
\newtheorem*{Remark 4}{Remark 4}

\numberwithin{equation}{section}

\begin{document}

\title[ Diffusion with a two-phase drift ]%
 {Transience, recurrence and speed of diffusions with a non-Markovian two-phase ``use it or lose it'' drift}

\author{ Ross G. Pinsky}
\address{Department of Mathematics\\
Technion---Israel Institute of Technology\\
Haifa, 32000\\ Israel}
\urladdr{http://www.math.technion.ac.il/~pinsky/}

\date{}

\begin{abstract}
We investigate the transience/recurrence  of a non-Markovian, one-dimensional diffusion process  which consists of a Brownian motion
with a non-anticipating drift that has two phases---a transient to $+\infty$  mode which is activated when the diffusion is sufficiently near its running maximum, and a recurrent mode
which is activated otherwise. We also consider the speed of  a diffusion with a two-phase drift, where the drift is equal to a certain positive constant when the diffusion
is sufficiently near its running maximum, and is equal to another positive constant otherwise.

\end{abstract}

\maketitle
\section{Introduction and Statement of Results}
Over the past fifteen years or so, there has been much interest in the study of the long term  behavior of various random walks with
non-Markovian transition mechanisms, such as random walks in random environment, self-avoiding random walk,
 edge or vertex reinforced random walks, and excited (so-called  ``cookie'') random walks.
See, for example, the monograph \cite{Z}, and the survey articles \cite{P07} and \cite{KZ12}, which include many references.
Non-Markovian diffusion processes analogous to excited random walks  have also  been studied
(see \cite{D96}, \cite{D99}, \cite{RS12}), as well as so-called Brownian polymers, which are non-Markovian self-repelling
diffusions, analogous to certain negatively reinforced random walks (see \cite{DR92}, \cite{CM96}, \cite{MT08}, \cite{HTT12}).
In this paper we investigate the transience/recurrence  of a non-Markovian, one-dimensional diffusion process  which consists of a Brownian motion
with a non-anticipating drift that has two phases---a transient to $+\infty$  mode which is activated when the diffusion is sufficiently near its running maximum, and a recurrent mode
which is activated otherwise.
We also consider the speed of  a diffusion with a two-phase drift, where the drift is equal to a certain positive constant when the diffusion
is sufficiently near its running maximum, and is equal to another positive constant otherwise.

Let $b^T(x)$ and $b^R(x)$ be continuous functions on $R$ which satisfy
\begin{equation}\label{1dt-r}
\begin{aligned}
&\int_{-\infty}\exp(-\int_0^x2b^T(y)dy)dx=\infty,\ \int^\infty\exp(-\int_0^x2b^T(y)dy)dx<\infty;\\
&\int_{-\infty}\exp(-\int_0^x2b^R(y)dy)dx=\infty,\ \int^\infty\exp(-\int_0^x2b^R(y)dy)dx=\infty.\
\end{aligned}
\end{equation}
As is well known \cite{RP95}, the one-dimensional diffusion processes  corresponding to the operators $\mathcal{L}^T\equiv\frac12\frac{d^2}{d x^2}+b^T(x)\frac d{dx}$
and
$\mathcal{L}^R\equiv\frac12\frac{d^2}{d x^2}+b^R(x)\frac d{dx}$
are respectively transient to $+\infty$ and recurrent.
Let $\gamma:[0,\infty)\to(0,\infty)$  be a continuous function  satisfying
$$
\gamma>0,\ \gamma'<1 \ \text{and}\ \lim_{x\to\infty}(x-\gamma(x))=\infty.
$$
For a continuous trajectory $x(\cdot):[0,\infty)\to R$, let $x^*(t)=\max_{0\le s\le t}x(s)$ denote its running maximum.
We define a non-anticipating, non-Markovian drift $b(t,x(\cdot))$ by
\begin{equation}\label{t-r}
b(t,x(\cdot))=\begin{cases} b^T(x(t)),\ \text{if}\ x(t)>x^*(t)-\gamma(x^*(t));\\
b^R(x(t)),\ \text{if}\ x(t)\le x^*(t)-\gamma(x^*(t)).\end{cases}
\end{equation}
We  consider the diffusion process $X(t)$  that satisfies the stochastic differential equation
\begin{equation}\label{sde}
X(t)=x_0+W(t)+\int_0^tb(s,X(\cdot))ds,
\end{equation}
where $W(\cdot)$ is a Brownian motion.
Existence and uniqueness for this stochastic differential equation follow from the standard theory for classical diffusion processes
(see section 3).
We call the solution to \eqref{sde}  a diffusion with a two-phase ``use it or lose it'' drift.

For example, the process $X(t)$ might represent the price of a stock.
As prices  rise, people are encouraged to buy, creating a certain trend represented by the transient drift.
In addition to this underlying trend, there is a random fluctuation represented by the Brownian motion.
These  random fluctuations might
cause  prices to  slump. If the slump becomes sufficiently large, it
discourages buying, which creates a new weaker trend, represented by the recurrent drift.
When random fluctuations eventually result in prices rising to levels  close to the previous high, the original stronger
trend reasserts itself.
With similar reasoning, the process $X(t)$ might measure   sales figures for a trendy brand product, or some other measurement
of the general zeitgeist.

We call $\gamma$ the ``down-crossing'' function.
For the majority of the paper, we will consider the case that the down-crossing function  $\gamma$ is a constant. In this case, at any time $t$, the diffusion $X(t)$
will run in the transient mode if and only if $X(t)> X^*(t)-\gamma$, or equivalently, if and only if by time $t$ the
path $X(\cdot)$  has not down-crossed an interval of length $\gamma$ whose left hand endpoint is larger than or equal to $X(t)$.

The various equivalent definitions of transience and recurrence for classical non-degenerate  diffusion processes
 hold for the diffusion with the two-phase drift. (This will be clear from the construction in section 3.)
We state here the standard definitions, although we will use other equivalent definitions
in the proofs. The diffusion with the two-phase drift is  \it recurrent\rm\ if for any pair of points $x_0$ and
$x_1$,  the process starting at $x_0$ almost surely returns to  $x_1$ at arbitrarily large times. The diffusion with the
two-phase drift is  \it transient to $+\infty$\rm\ if
starting at any $x_0$, the process almost surely  satisfies $\lim_{t\to\infty}X(t)=\infty$.

Our first result concerns the case in which the transient drift is constant: $b^T(x)\equiv b>0$, and the recurrent drift
$b^R$ satisfies a regularity condition; namely, that the drift $b^R\vee 0$ is also a recurrent drift. In this case, the diffusion with the two-phase drift
is always recurrent.
\begin{theorem}\label{alwaysrec}
Assume that the down-crossing function $\gamma$ is constant.
Let the transient drift be constant: $b^T(x)\equiv b>0$, and assume that the recurrent drift $b^R$ is such that
the drift $b^R\vee0$ is also recurrent. That is, assume that the condition satisfied by $b^R$ in \eqref{t-r} is also satisfied
by $b^R\vee0$. Then the diffusion with the two-phase drift is recurrent.
\end{theorem}

\noindent\bf Remark.\rm\ Note that the diffusion with the two-phase drift is recurrent
even if the recurrent drift $b^R$ is just  border line recurrent---for example, if for sufficiently large $x$, $b^R(x)=\frac1{2x}$ is
 the drift of the radial part of a two-dimensional Brownian motion.

Maintaining the constant transient drift, but
 choosing the  recurrent drift $b^R$ to take on very large positive values at most locations, and compensating to insure recurrence by having
it take on even much larger   negative values at  other locations, we can  construct a  diffusion with such a two-phase drift
that is transient.

\begin{theorem}\label{exampletrans}
Assume that the down-crossing function $\gamma$ is constant.
Let the transient drift be constant: $b^T(x)=b>0$.
There exists a  recurrent drift $b^R$ such that the corresponding diffusion with the two-phase drift
is transient.
\end{theorem}

We continue to assume that the down-crossing function $\gamma$ is constant.
As noted above,  $b(t,X(\cdot))=b^T(X(t))$ if and only if by time $t$, the
path $X(\cdot)$  has not down-crossed an interval of length $\gamma$ whose left hand endpoint is larger than or equal to $X(t)$.
Now if
the transient diffusion corresponding to $\mathcal{L}^T$ is such that it almost surely eventually stops making down-crossings of length $\gamma$, then
the diffusion in the two-phase environment will eventually stop making down-crossings of length $\gamma$  and will eventually
be driven just by the transient drift; consequently, it will be transient. In \cite{RP10} we proved the following result.
\medskip

\noindent \bf Theorem P2.\it\
Consider the diffusion process  corresponding to
the operator $\mathcal{L}^T$.

\noindent i. If $b^T(x)\le \frac1{2\gamma}\log x+\frac1{\gamma}\log^{(2)} x$ for sufficiently large $x$, then the diffusion  almost
surely makes
$\gamma$-down-crossings for arbitrarily large times;

\noindent ii. If $b^T(x)\ge\frac1{2\gamma}\log x+\frac k\gamma\log^{(2)} x$, for some $k>1$
and for sufficiently large $x$, then the diffusion almost surely
 eventually stops making down-crossings of size $\gamma$.
\rm \medskip

In light of Theorem P2 and the paragraph preceding it, in the case of a constant down-crossing function $\gamma$, if $b^T$ satisfies the condition in part (ii) of the theorem, then the diffusion
with the two-phase drift is transient, regardless of what the recurrent drift $b^R$ is.

Continuing with a constant down-crossing function $\gamma$, we now let the recurrent diffusion be Brownian motion, that is, $b^R\equiv0$, and determine what the threshold growth rate is on $b^T$
that distinguishes between transience and recurrence for the diffusion with the two-phase drift.
\begin{theorem}\label{driftdisting}
Assume that the down-crossing function $\gamma$ is constant. Let the recurrent diffusion be Brownian motion: $b^R\equiv0$.

\noindent i. If $b^T(x)\le \frac1{2\gamma}\log^{(2)}x+\frac1{2\gamma}\log^{(3)}x, \ \text{for large}\ x$, then the diffusion with the two-phase drift
is recurrent;

\noindent ii. If $b^T(x)\ge \frac1{2\gamma}\log^{(2)}x+\frac k{2\gamma}\log^{(3)} x,\ \text{for large}\ x$, where $k>1$,  then the
diffusion with the two-phase drift
is transient.
\end{theorem}

We now consider the case that the recurrent diffusion is Brownian motion: $b^R\equiv0$, that the transient drift is constant: $b^T(x)\equiv b$,
but we allow the down-crossing function $\gamma$ to grow with $x$. We determine the threshold growth rate on the down-crossing function
that distinguishes between transience and recurrence for the diffusion with the two-phase drift.
\begin{theorem}\label{down-crossingdisting}
Let the recurrent diffusion be Brownian motion: $b^R\equiv0$, and let the transient drift be constant: $b^T(x)\equiv b$.

\noindent i. If the down-crossing function $\gamma$ satisfies
 $\gamma(x)\le \frac1{2b}\log^{(2)}x+\frac1{2b}\log^{(3)}x, \ \text{for large}\ x$, then the diffusion with the two-phase drift
is recurrent;

\noindent ii. If the down-crossing function $\gamma$ satisfies  $\gamma(x)\ge \frac1{2b}\log^{(2)}x+\frac k{2b}\log^{(3)} x,\ \text{for large}\ x$, where $k>1$,  then the
diffusion with the two-phase drift
is transient.

\end{theorem}
\medskip

We now  make one major and one minor change in the setup we have used until now. Let the two-phase drift $b(t,x(\cdot))$ be given by
\eqref{t-r}, with $b^T\equiv b$, with $\gamma$ constant, and with the recurrent drift $b^R$ replaced by
the transient constant drift $c$. The case $c\in(0,b)$ is more in keeping with the theme of this paper, but
the case $c>b$ also has an interesting aspect. In particular, for the case $c=\infty$, described below,  the dependence on the diffusion
coefficient is worth noting; thus, we will consider the operator
\begin{equation}\label{ballisticop}
\begin{aligned}
&L=\frac12a\frac{d^2}{dx^2}+b(t,X(\cdot))\frac d{dx},\\
&\text{where}\ b(t,X(\cdot))=\begin{cases} b>0,\ \text{if}\ X(t)>X^*(t)-\gamma;\\ c>0,\ \text{if}\ X(t)\le X^*(t)-\gamma.
\end{cases}
\end{aligned}
\end{equation}
The next theorem gives the speed of the diffusion in this two-phase drift.
 \begin{theorem}\label{ballistic}
Consider the diffusion in the two-phase drift corresponding to the operator $L$ in \eqref{ballisticop}.
The speed of the diffusion $X(t)$ with the two-phase drift  is given by
$$
\lim_{t\to\infty}\frac{X(t)}t=\Big(\frac{c(e^{\frac{2 b\gamma}a}-1)}{c(e^{\frac{2 b\gamma}a}-1)+( b-c)(1-e^{-\frac{2 b\gamma}a})}\Big)b\ \ \text{a.s.}
$$
 \end{theorem}
\noindent \bf Remark.\rm\ Let  $d(b,c,\gamma,a)\equiv\frac{c(e^{\frac {2b\gamma}a}-1)}{c(e^{\frac{2 b\gamma}a}-1)+( b-c)(1-e^{-\frac{2 b\gamma}a})}$.
In the case $c\in(0,b)$, we call $d(b,c,\gamma,a)$ the  \it damping coefficient\rm. It gives the fractional reduction in speed between a classical diffusion
with drift $b$ and the slowed down diffusion with two-phase drift---$b$ when the process is less than
 distance $\gamma$ from its running maximum, and  $c$   when the process is
at least  distance $\gamma$ from its running maximum. Of course, it is clear that $d(b,c,\gamma)$ must always fall between
$\frac cb$ and 1 when $c\in(0,b)$.
We make the following observations:

\noindent 1. When $b\to\infty$, the damping coefficient $d(b,c,\gamma,a)$  converges to 1 exponentially  in $b$;

\noindent 2. When $\gamma\to\infty$,  $d(b,c,\gamma,a)$ converges to 1 exponentially in $\gamma$;

\noindent 3. When $c\to0$, the damping coefficient $d(b,c,\gamma,a)$ converges to 0 linearly in $c$;

\noindent 4. When $c\to0$ and $\gamma\to\infty$, the damping coefficient $d(b,c,\gamma,a)$ converges to 1 (respectively, converges to 0, remains bounded away
from 0 and 1) if $c\exp(\frac{2 b\gamma}a)$ converges to $\infty$ (respectively, converges to 0, remains bounded away from 0 and $\infty$);

\noindent 5. When $c\to0$ and $b\to\infty$, the damping coefficient $d(b,c,\gamma,a)$ converges to 0 if $c\exp(\frac{2 b\gamma}a)$ remains bounded.
Otherwise, the damping coefficient converges to 0 (respectively, converges to 1, remains bounded away from 0 and 1)
if $\frac{c\exp(\frac{2 b\gamma}a)}b$ converges to 0 (respectively, converges to $\infty$, remains bounded away from 0 and $\infty$).

\noindent 6.  The limiting case $c=\infty$ corresponds to the situation in which
the value $X^*(t)-\gamma$ serves as a reflecting barrier for $X(t)$; thus, the process can never get farther than a distance $\gamma$ from
its running maximum.
We have $d(b,\infty,\gamma,a)=\frac{(e^{\frac{2 b\gamma}a}-1)}{(e^{\frac{2 b\gamma}a}+e^{-\frac{2 b\gamma}a}-2)}$.
For this process, as $\gamma\to\infty$, the speed satisfies  $d(b,\infty,\gamma,a)b\sim\frac a{2\gamma}$.
In particular, the speed approaches infinity and the leading order term is independent of $b$---the process is propelled
forward by the positive excursions of the Brownian motion with diffusion coefficient $a$.

\medskip

In section 2, we present some preliminary information on down-crossings. In section 3, we give an explicit
representation of the diffusion with a two-phase drift in terms of classical diffusions.
In section 4 we give a workable analytic criterion for transience/recurrence which depends on an auxiliary
discrete time, increasing Markov process. Sections 5-9 give the proofs of Theorems \ref{alwaysrec}-\ref{ballistic}.

\section{Preliminaries concerning down-crossings}
From the definition of the down-crossing function $\gamma$, it follows that  $x-\gamma(x)$ is increasing.
Define the stopping time  $\sigma_\gamma$ on
continuous paths $x(\cdot):[0,\infty)\to R$ with the standard filtration $\mathcal{F}_t=\sigma(x(s):0\le s\le t)$
by
$$
\begin{aligned}
&\sigma_\gamma=\inf\{t\ge0: \exists s<t \ \text{with}\ x(t)\le x(s)-\gamma(x(s))\}=\\
&\inf\{t\ge0:  x(t)= x^*(t)-\gamma(x^*(t))\}.
\end{aligned}
$$
The equality above follows from the fact that $x-\gamma(x)$ is increasing.
Let
$$
\begin{aligned}
&L^\gamma=x^*(\sigma_\gamma);\\
&K^\gamma=x(\sigma_\gamma)=x^*(\sigma_\gamma)-\gamma(x^*(\sigma_\gamma)).
\end{aligned}
$$

In the case that the down-crossing function $\gamma$ is constant,
$\sigma_\gamma$  is the first time that the path $x(\cdot)$ completes a down-crossing of an interval of length $\gamma$.
The interval that was down-crossed is $[K^\gamma,L^\gamma]$.
 In \cite{RP10}, for $\gamma$  constant,
  $L^\gamma$ was called  the \it $\gamma$-down-crossed
onset location.\rm\ In this paper, we will  use this terminology also for
variable $\gamma$. We will call $\sigma_\gamma$ the \it first $\gamma$-down-crossed time.\rm\
For use a bit later, let $\tau_a=\inf\{t\ge0: x(t)=a\}$ denote the first hitting time of the point $a$.

Consider now the one-dimensional diffusion process $Y(t)$ which corresponds to the operator $\mathcal{L}^T$ and which is transient to $+\infty$.
Denote probabilities for the process starting from $x$ by $P_x$.
Fixing a point $z_0$, let
\begin{equation}\label{uTfunc}
u_T(x)=\int_{z_0}^x\exp(-\int_{z_0}^y2b^T(r)dr)dy.
\end{equation}
(The formula in the theorem below is independent of $z_0$, but this specification of $z_0$ will be useful later on.)
The following result  was proved in \cite{RP10}.
\medskip

\noindent \bf Theorem P1.\it\  For the diffusion process corresponding to $\mathcal{L}^T$, and for constant $\gamma$,
the distribution of the $\gamma$-down-crossed onset location $L^\gamma$ is given by
$$
P_x(L^\gamma>x+y)=\exp(-\int_x^{x+y}\frac{u_T'(z)}{u_T(z)-u_T(z-\gamma)}dz), \ y>0.
$$
\medskip

\bf\noindent Remark 1.\rm\ In \cite{RP10}, where the notation is a bit different from here,
the mathematical definition of  $\sigma_\gamma$, and consequently also of $L^\gamma$, were written  incorrectly.
From the verbal description in \cite{RP10}, it is clear that the intended definition of $L^\gamma$ is the one
given here. All the proofs in \cite{RP10} are based on the correct definitions given here.

\bf\noindent Remark 2.\rm\ Theorem P2 in section 1 was  proved in \cite{RP10} as an application of Theorem P1.

The  same method of proof used to prove Theorem P1 can be used in the case of variable $\gamma$ to obtain a corresponding formula
for the distribution of $L^\gamma$.
\begin{proposition}\label{Lgamma}
For the diffusion process corresponding to $\mathcal{L}^T$, and for variable $\gamma$,
the distribution of the $\gamma$-down-crossed onset location $L^\gamma$ is given by
$$
P_x(L^\gamma>x+y)=\exp(-\int_x^{x+y}\frac{u_T'(z)}{u_T(z)-u_T(z-\gamma(z))}dz), \ y>0.
$$
\end{proposition}
\medskip

\section{A representation for  diffusion with  two-phase drift}

Consider the diffusion $X(t)$ with the two-phase drift starting from $x_0$.
Up until the first $\gamma$-down-crossed time  $\sigma_\gamma$, the process is exactly the $Y(\cdot)$-process corresponding to the operator $\mathcal{L}^T$
and starting from $x_0$.
We have $X(\sigma_\gamma)=Y(\sigma_\gamma)=K^\gamma$
and $X^*(\sigma_\gamma)=Y^*(\sigma_\gamma)=L^\gamma$, with
$L^\gamma$ distributed as in Proposition \ref{Lgamma}.
Let
$$
\hat\tau_{L^\gamma}=\inf\{t\ge0: X(\sigma_\gamma+t)=L^\gamma\}.
$$
Then
$\sigma_\gamma+\hat\tau_{L^\gamma}$
is the first time after $\sigma_\gamma$ that the process $X(\cdot)$  returns to its running maximum $L^\gamma$.
Let $Z^{R,z,T}(t)$ be the diffusion starting from $z$ and corresponding to the operator
$\mathcal{L}^{R,z,T}=\frac12\frac{d^2}{dx^2}+b^{R,z,T}(x)\frac d{dx}$, where
$$
b^{R,z,T}(x)=\begin{cases} b^R(x),\ x\le z;\\ b^T(x),\ x>z.\end{cases}
$$
Then the distribution of $\{X(\sigma_\gamma+t), 0\le t\le \hat\tau_{L^\gamma}\}$, conditioned on
 $K^\gamma=X(\sigma_\gamma)=z$
and $L^\gamma=X^*(\sigma_\gamma)=a$,
is that of $\{Z^{R,z,T}(t), 0\le t\le \tau_a\}$.
Of course, $X(\sigma_\gamma+\hat\tau_{L^\gamma})=X^*(\sigma_\gamma+\hat\tau_{L^\gamma})=L^\gamma$.
Starting from time $\sigma_\gamma+\hat\tau_{L^\gamma}$, when the process has returned to its running maximum,  $X$ again looks like the
$Y$ process, until it again performs a $\gamma$-down-crossing, at which point it becomes
a $Z^{R,z,T}$ process for appropriate $z$ until it returns to its running maximum, and everything is repeated again.

In light of the above description, $X(t)$  can be  represented as follows.
For each $x\in R$ and each $n\ge1$,  let $Y^{n,x}(\cdot)$ be a diffusion process corresponding to the operator
$\mathcal{L}^T$ and starting from $x$. Make the processes independent for different pairs $(n,x)$.
Similarly, for each $z\in R$ and each $n\ge1$, let
$Z^{R,z,T,n}(\cdot)$ be a diffusion process corresponding to the operator
$\mathcal{L}^{R,z,T}$ and starting from $z$. Make the processes independent for different pairs $(n,z)$ and independent of the $Y^{n,x}$ processes.
Let $\sigma^{n,x}_\gamma$ denote the first $\gamma$-down-crossing time for the process $Y^{n,x}$, and let $\tau_a^{n,z}$ denote
the first hitting time of $a$ for the process $Z^{R,z,T,n}$.
Let $L_0^\gamma=x_0$ and
 then by induction, for $n\ge1$, define
$L_n^\gamma$ to be the $\gamma$-down-crossed onset location  for $Y^{n,L_{n-1}^\gamma}$.
For $n\ge1$, let $K_n^\gamma$ correspond to  $L_n^\gamma$ as $K^\gamma$ corresponds to  $L^\gamma$.
Then $X(\cdot)$ can be represented as
\begin{equation}\label{represent1}
\begin{aligned}
&X(t)=Y^{1,x_0}(t),\ 0\le t\le \sigma_\gamma^{1,x_0}~;\\
&X(t)=Z^{R,K_1^\gamma,T,1}(t),\ \sigma_\gamma^{1,x_0}\le t\le
\sigma_\gamma^{1,x_0}+\tau^{1,K_1^\gamma}_{L_1^\gamma},
\end{aligned}
\end{equation}
and for $n\ge2$,
\begin{equation}\label{represent2}
\begin{aligned}
&X(\sum_{j=1}^{n-1}\sigma_\gamma^{j,L_{j-1}^\gamma}+\sum_{j=1}^{n-1}\tau_{L_j^\gamma}^{j,K_j^\gamma}+t)=Y^{n,L_{n-1}^\gamma}(t),\ 0\le t\le
\sigma_\gamma^{n,L_{n-1}^\gamma};\\
&X(\sum_{j=1}^n\sigma_\gamma^{j,L_{j-1}^\gamma}+\sum_{j=1}^{n-1}\tau_{L_j^\gamma}^{j,K_j^\gamma}+t)=Z^{R,K_n^\gamma,T,n}(t),\ 0
\le t\le\tau_{L_n^\gamma}^{n,K_n^\gamma}.
\end{aligned}
\end{equation}
From the above representation, it is clear that existence and uniqueness for \eqref{sde} follows from the standard theory for
classical diffusion processes.

At those times $s$ when $X(s)$ is running as a $Y^{n,L_{n-1}^\gamma}(\cdot)$-process,  for some $n\ge1$, we will say that $X(\cdot)$ is in the $Y$-mode,
and at those times $s$ when $X(s)$ is running as a $Z^{R,K_n^\gamma,T,n}(\cdot)$-process, for some $n\ge1$, we will say
that $X(\cdot)$ is in the $Z$-mode.
We denote by  $\mathcal{P}_{x_0}$  probabilities corresponding to the diffusion $X(t)$ with
the two-phase drift starting from $x_0$, and by $\mathcal{E}_{x_0}$ the corresponding expectation.

Note that  $\{L_n^\gamma\}_{n=0}^\infty$ is a monotone increasing Markov process
whose  transition distribution from a state $x$
is the distribution of $L^\gamma$ in Proposition \ref{Lgamma}. That is,
the transition probability measure $p(x,\cdot)$ is given by
\begin{equation}\label{transpd}
\begin{aligned}
&p(x,[x+y,\infty))=
\exp(-\int_x^{x+y}\frac{u_T'(z)}{u_T(z)-u_T(z-\gamma(z))}dz),
\text{for}\ y\ge0.
\end{aligned}
\end{equation}
We will use the notation $P^L_{x_0}$ to denote probabilities for this discrete time Markov process and will denote the corresponding expectation
by $E_{x_0}^L$.
The  times  $\{L_n^\gamma\}_{n=0}^\infty$ will be called   ``regeneration'' points for $X(\cdot)$  because
the $\mathcal{P}_{x_0}$-distribution of $\{X(\sum_{j=1}^{n-1}\sigma_\gamma^{j,L_{j-1}^\gamma}+\sum_{j=1}^{n-1}\tau_{L_j^\gamma}^{j,K_j^\gamma}+t), 0\le t<\infty\}$ given that $L_{n-1}^\gamma=a$   is the same as the
$\mathcal{P}_a$-distribution of \{$X(t), 0\le t<\infty\}$.

\section{A transience/recurrence criterion}
From the construction in section 3,
the well-known equivalent alternative conditions for transience/recurrence, which hold for standard non-degenerate diffusion
processes \cite{RP95}, are easily seen to hold for diffusions with a two-phase drift. Since it is clear that
the process is either transient to $+\infty$ or recurrent,
we can use the following criterion:

\begin{equation}\label{t-r-crit}
\begin{aligned}
&\text{{\bf Transience}: for some pair of points}\ z_0<x_0,\ \text{one has}\ \mathcal{P}_{x_0}(\tau_{z_0}=\infty)>0;\\
&\text{{\bf Recurrence}: for some pair of points}\ z_0<x_0,  \ \text{one has}\ \mathcal{P}_{x_0}(\tau_{z_0}=\infty)=0.
\end{aligned}
\end{equation}

We choose the point $z_0$ so that $z_0<x_0-\gamma(x_0)$. Then it follows from the representation  of the process $X(\cdot)$
 that at time $\tau_{z_0}$ the process is in the $Z$-mode.
Using this with  the regeneration structure noted at the end of the previous section, it follows that
\begin{equation}\label{t-r-explicit}
\mathcal{P}_{x_0}(\tau_{z_0}=\infty)=E_{x_0}^LG(\{L_n^\gamma\}_{n=1}^\infty),
\end{equation}
where for any nondecreasing sequence $\{a_n\}_{n=1}^\infty$ satisfying
$a_1\ge x_0$, we define
\begin{equation}\label{G-func}
G(\{a_n\}_{n=1}^\infty)=\prod_{n=1}^\infty P_{a_n-\gamma(a_n)}^{R,a_n-\gamma(a_n),T}(\tau_{a_n}<\tau_{z_0}),
\end{equation}
and where $P_z^{R,z,T}$ denotes probabilities for a $Z^{R,z,T}$-processes starting at $z$.
From \eqref{t-r-explicit} we conclude  that $\mathcal{P}_{x_0}(\tau_{z_0}=\infty)=0$ if and only if
$G(\{L_n^\gamma\}_{n=1}^\infty)=0$ a.s. $P_{x_0}^L$; thus, from \eqref{G-func} we have
\begin{equation}\label{recurr-prod-cond}
\mathcal{P}_{x_0}(\tau_{z_0}=\infty)=0\ \
\text{if and only if}\
\sum_{n=1}^\infty\Big(P_{L_n^\gamma-\gamma(L_n^\gamma)}^{R,L_n^\gamma-\gamma(L_n^\gamma),T}(\tau_{z_0}<\tau_{L_n^\gamma})\Big)=\infty\
\text{a.s.}\ P_{x_0}^L.
\end{equation}

As is well known \cite{RP95}, the function $v(x)\equiv P_x^{R,z,T}(\tau_{z_0}<\tau_{z+c})$, for $z_0\le x\le z+c$,
satisfies $\mathcal{L}^{R,z,T}v=0$ in $(z_0,z)\cup(z,z+c)$, $v(z_0)=1, v(z+c)=0$, $v(z^-)=v(z^+)$, and $v'(z^-)=v'(z^+)$.
Solving this, one finds that $P_z^{R,z,T}(\tau_{z_0}<\tau_{z+c})=v(z)$ is given by
\begin{equation}\label{hittingtimeprob}
\begin{aligned}
&P_z^{R,z,T}(\tau_{z_0}<\tau_{z+c})=\\
&\frac{\exp(-\int_{z_0}^z(2b^R)(y)dy)\int_z^{z+c}dy\exp(-\int_z^y2b^T(r)dr)}{\int_{z_0}^zdy \exp(-\int_{z_0}^y2b^R(r)dr)+
\exp(-\int_{z_0}^z(2b^R)(y)dy)\int_z^{z+c}dy\exp(-\int_z^y2b^T(r)dr)}.
\end{aligned}
\end{equation}
Analogous to $u_T$ in \eqref{uTfunc}, define  the function
\begin{equation}\label{uRfunc}
u_R(x)=\int_{z_0}^x\exp(-\int_{z_0}^y2b^R(r)dr)dy.
\end{equation}
We  then rewrite the somewhat unwieldy equation \eqref{hittingtimeprob}, which would become a lot more unwieldy below, in the form
\begin{equation}\label{wieldy}
P_z^{R,z,T}(\tau_{z_0}<\tau_{z+c})=\frac{u'_R(z)\big(u_T(z+c)-u_T(z)\big)
\exp(\int_{z_0}^z2b^T(y)dy)}{u_R(z)+u'_R(z)\big(u_T(z+c)-u_T(z)\big)\exp(\int_{z_0}^z2b^T(y)dy)}.
\end{equation}

From \eqref{t-r-crit}, \eqref{recurr-prod-cond} and \eqref{wieldy}, we obtain the following transience/recurrence criterion for the diffusion
with the two-phase drift.
\begin{proposition}\label{applicable-t-r}
Define
\begin{equation}\label{applicablecond}
\begin{aligned}
&H(s)=\\
&\frac{u'_R(s-\gamma(s))\big(u_T(s)-u_T(s-\gamma(s))\big)\exp(\int_{z_0}^{s-\gamma(s)}
2b^T(y)dy)}{u_R(s-\gamma(s))+u'_R(s-\gamma(s))\big(u_T(s)-u_T(s-\gamma(s))\big)
\exp(\int_{z_0}^{s-\gamma(s)}2b^T(y)dy)}.
\end{aligned}
\end{equation}
Let $\{L_n^\gamma\}_{n=0}^\infty$ be the monotone increasing Markov process with $L_0^\gamma=x_0$ and transition probability measure
$p(x,\cdot)$ given by \eqref{transpd}. If
\begin{equation}
 \sum_{n=0}^\infty H(L^\gamma_n)=
\infty \ \text{a.s.},
\end{equation}
then the diffusion with the  two-phase drift is recurrent. Otherwise, it is transient.
\end{proposition}
\section{Proof of Theorem \ref{alwaysrec}}
\it \noindent Proof of Theorem \ref{alwaysrec}.\rm\
By the Ikeda-Watanabe  comparison theorem \cite{IW89}, if we prove recurrence for the diffusion with the two-phase drift
whose recurrent drift $b^R$ is replaced by the drift $b^R\vee0$, then original diffusion with the two-phase drift
is also recurrent. By assumption, the drift $b^R\vee0$ is also a recurrent drift. Thus, we may assume without loss
of generality that the recurrent drift $b^R$ is nonnegative.

Since $b^T\equiv b$, we have from \eqref{uTfunc} that
$$
u_T(x)=\frac1{2b}\Big(1-\exp\big(-2b(x-z_0)\big)\Big).
$$
Using this along with the fact that $\gamma$ is constant, we have
\begin{equation}\label{uTbconst}
\big(u_T(s)-u_T(s-\gamma(s))\big)\exp(\int_{z_0}^{s-\gamma(s)}2b^T(y)dy)=\frac1{2b}(1-\exp(-2b\gamma))\equiv c_{b,\gamma},
\end{equation}
 and thus
the formula for $H(s)$ in \eqref{applicablecond} simplifies
to
\begin{equation}\label{H1}
H(s)=\frac{c_{b,\gamma}u'_R(s-\gamma)}{u_R(s-\gamma)+c_{b,\gamma}u'_R(s-\gamma)}.
\end{equation}
From \eqref{uRfunc}, it follows that $u_R$ is increasing, and by the assumption that $b^R$ is nonnegative, it follows that
$u'_R$ is non-increasing.  Thus, $H$ is non-increasing.

We also have
\begin{equation}\label{dbgamma}
\frac{u_T'(z)}{u_T(z)-u_T(z-\gamma(z))}=\frac1{c_{b,\gamma}\exp(2b\gamma)}\equiv\frac1{ d_{b,\gamma}},
\end{equation}
and thus the increment measure $p(x,x+A)$, $A\subset[0,\infty)$, corresponding to the
 transition probability measure $p(x,\cdot)$ in \eqref{transpd} for the Markov process $\{L_n^\gamma\}_{n=0}^\infty$ is
independent of $x$ and is equal to the exponential density with parameter $\frac1{d_{b,\gamma}}$. Consequently,
$L_n^\gamma-L_0^\gamma$ is the sum of $n$ IID exponential random variables with the above parameter, and thus
\begin{equation}\label{LLN}
\lim_{n\to\infty}\frac{L^\gamma_n}n=d_{b,\gamma} \ \text{a.s.}
\end{equation}
Using \eqref{LLN} along with the fact that $H$ is non-increasing,
if we show that $\sum_{n=1}^\infty H((d_{b,\gamma}+1)n)=\infty$, then
it follows that
$\sum_{n=1}^\infty H(L^\gamma_n)=
\infty \ \text{a.s.}$, and consequently, from Proposition \ref{applicable-t-r} we conclude that the diffusion with the two-phase
drift is recurrent.
Since $u_R'$ is non-increasing and $u_R$ is increasing, it is easy to see from \eqref{H1} that
$\sum_{n=1}^\infty H((d_{b,\gamma}+1)n)=\infty$ if and only if $\sum_{n=1}^\infty \hat H((d_{b,\gamma}+1)n)=\infty$,
where $\hat H=\frac{u'_R(s-\gamma)}{u_R(s-\gamma)}$. Since $\hat H$ is monotone, it follows that
$\sum_{n=1}^\infty \hat H((d_{b,\gamma}+1)n)=\infty$ if and only if $\int^\infty\frac{u'_R(s)}{u_R(s)}ds=\infty$,
that is, if and only if $\lim_{s\to\infty}u_R(s)=\infty$. But this last inequality holds from \eqref{uRfunc} and \eqref{1dt-r} since $b^R$ is a recurrent drift.
\hfill $\square$

\section{Proof of Theorem \ref{exampletrans}}
As in the proof of Theorem \ref{alwaysrec}, the random variables $\{L_n^\gamma-L_{n-1}^\gamma\}_{n=1}^\infty$ are IID and distributed
according to the exponential distribution with parameter $\frac1{d_{b,\gamma}}$, and the function $H$ is given by \eqref{H1}.
In order to show transience, by Proposition \ref{applicable-t-r}, we need to show that
$\sum_{n=0}^\infty H(L_n^\gamma)<\infty$ with positive probability.
Recalling \eqref{uRfunc} and \eqref{1dt-r}, to complete the proof, we
 will construct a function $u_R$ which satisfies $u_R>0, u_R'>0$ and $\lim_{x\to\infty}u_R(x)=\infty$, and for which the above sum
is almost surely finite.
(The corresponding drift $b^R$ will then be given by $-\frac{u_R''}{u_R}$.)

For $j\ge2$, define the interval $I_j=[x_0-\gamma+j,x_0-\gamma+ j+\frac1{j^2}]$.
Without loss of generality, assume that $x_0-\gamma+2\ge0$.
We now show that
\begin{equation}\label{B-C}
P_{x_0}^L(L_n^\gamma-\gamma\in\cup_{j=2}^\infty I_j)\le \frac c{n^2},
\end{equation}
for some $c>0$.
The distribution of $L_n^\gamma-x_0$ is that of the sum of $n$ IID exponential random variables with parameter
$\lambda\equiv\frac1{d_{b,\gamma}}$.
Thus, its density function is $\frac{\lambda^nx^{n-1}}{(n-1)!}\exp(-\lambda x)$, $x\ge0$.
For an appropriate constant $C>1$, we then have for $n\ge 3$,
\begin{equation}\label{intervals}
\begin{aligned}
&P_{x_0}^L(L_n^\gamma-\gamma\in\cup_{j=2}^\infty I_j)=\sum_{j=2}^\infty\int_{I_j}\frac{\lambda^nx^{n-1}}{(n-1)!}\exp(-\lambda x)dx\le\\
&C\int_0^\infty\frac{\lambda^nx^{n-3}}{(n-1)!}\exp(-\lambda x)dx=\frac{C\lambda^2}{(n-2)(n-2)}.
\end{aligned}
\end{equation}
Now \eqref{B-C} follows from \eqref{intervals}.

We now construct a positive, strictly increasing $C^1$-function $u_R$ whose derivative on $R-\cup_{j=2}^\infty I_j$ is uniformly bounded,
and which satisfies $u_R(x)\ge x^2$, for $x\ge2$. (Of course, to have this quadratic growth, $u'_R$ must get very large
at certain places on $\cup_{j=2}^\infty I_j$.)
By the law of large numbers, $\lim_{n\to\infty}\frac{L_n^\gamma}n=d_{b,\gamma}$ a.s.
By \eqref{B-C} and  the lemma of Borel-Cantelli, $P_{x_0}^L(L_n^\gamma-\gamma\in\cup_{j=2}^\infty I_j\ \text{i.o.})=0$.
Using the facts noted in this paragraph, we conclude
that
$$
\sum_{n=0}^\infty H(L_n^\gamma)=\sum_{n=0}^\infty\frac{c_{b,\gamma}u'_R(L_n^\gamma-\gamma)}{u_R(L_n^\gamma-\gamma)+
c_{b,\gamma}u'_R(L_n^\gamma-\gamma)}<\infty\ \text{a.s.}
$$

\hfill $\square$

\section{Proof of Theorem \ref{driftdisting}}
Since $b^R=0$, we have from \eqref{uRfunc} that $u^R(x)=x-z_0$. Also, $\gamma$
is constant.
Thus, from \eqref{applicablecond}, we have
\begin{equation}\label{H-driftdisting}
H(s)=\frac{\int_{s-\gamma}^sdy\exp(-\int_{s-\gamma}^y2b^T(r)dr)}{s-\gamma+\int_{s-\gamma}^sdy\exp(-\int_{s-\gamma}^y2b^T(r)dr)}.
\end{equation}

By comparison, it suffices to consider the case that
$$
b^T(x)=\frac1{2\gamma}\log^{(2)}x+\frac k{2\gamma}\log^{(3)}x,
$$
 for $x\ge z_0$, with $z_0$ large enough so that $\log^{(3)}z_0$ is defined.
We need to show transience for $k>1$ and recurrence for $k=1$.
In what follows, we will always assume that $k\ge1$.

We have
\begin{equation}\label{logint}
(y-s+\gamma\big)\log^{(j)}(s-\gamma)\le \int_{s-\gamma}^y\log^{(j)}rdr\le (y-s+\gamma\big)\log^{(j)}s,\ s-\gamma\le y\le s.
\end{equation}
Thus,
\begin{equation}\label{est-s-g-s}
\begin{aligned}
&\frac{\gamma-o(s)}{\log^{(2)}s+k\log^{(3)}s}\le \int_{s-\gamma}^sdy\exp(-\int_{s-\gamma}^y2b^T(r)dr)\le \\ &\frac\gamma{\log^{(2)}(s-\gamma)+k\log^{(3)}(s-\gamma)},\ \text{as}\ s\to\infty.
\end{aligned}
\end{equation}
From \eqref{H-driftdisting} and \eqref{est-s-g-s} we conclude that there exist constants $C_1,C_2>0$ such that
\begin{equation}\label{H-est-driftdisting}
\frac{C_1}{s\log^{(2)}s}\le H(s)\le \frac{C_2}{s\log^{(2)}s}, \ \text{for large}\ s.
\end{equation}

We now investigate the growth rate of the Markov process $\{L_n^\gamma\}_{n=0}^\infty$. Recall that given $L_j^\gamma=x$, the distribution
of $L_{j+1}^\gamma-L_j^\gamma$
 is the distribution given in \eqref{transpd}.
From \eqref{uTfunc} we have
\begin{equation}\label{altform}
\begin{aligned}
&\frac{u_T'(x)}{u_T(x)-u_T(x-\gamma)}=\frac{\exp(-\int_{z_0}^x2b^T(r)dr)}{\int_{x-\gamma}^xdy\exp(-\int_{z_0}^y2b^T(r)dr)}=\\
&\frac{\exp(-\int_{x-\gamma}^x2b^T(r)dr)}{\int_{x-\gamma}^xdy\exp(-\int_{x-\gamma}^y2b^T(r)dr)}.
\end{aligned}
\end{equation}
From \eqref{logint} and the definition of $b^T$, we have
$$
\frac{c_1}{(\log x)(\log^{(2)}x)^k}\le\exp(-\int_{x-\gamma}^x2b^T(r)dr)\le \frac{c_2}{(\log x)(\log^{(2)}x)^k}, \ \text{for large}\ x,
$$
for constants $c_1,c_2>0$. Using this with
\eqref{est-s-g-s} and \eqref{altform}, we conclude that there exist constants $C_3,C_4>0$ such that
\begin{equation}\label{expparameter}
 \frac{C_3}{(\log x)(\log^{(2)}x)^{k-1}}\le \frac{u_T'(x)}{u_T(x)-u_T(x-\gamma)}\le \frac{C_4}{(\log x)(\log^{(2)}x)^{k-1}}, \ \text{for large}\ x.
\end{equation}

We now prove recurrence in the case that $k=1$.
Let $\{x_n\}_{n=0}^\infty$ be a sequence of positive numbers with $x_0$ sufficiently large so that
the bound in \eqref{expparameter} holds for $x\ge x_0$, and let $s_n=\sum_{j=0}^nx_j$.
We have $\int_{s_{n-1}}^{s_n} \frac1{\log z}dz\ge \frac {x_n}{\log s_n}$.
Using this with \eqref{expparameter} and \eqref{transpd}, we have
\begin{equation}\label{cond-L}
P_{x_0}^L(L^\gamma_n-L_{n-1}^\gamma\le  x_n|L^\gamma_{n-1}\le s_{n-1})\ge1-\exp\big(-\frac {C_3x_n}{\log s_n}\big),\ n\ge1.
\end{equation}
Fix  a large number $M$.
With $x_0$ as above, we  wish to select the sequence $\{x_n\}_{n=0}^\infty$ so that
\begin{equation}\label{generatingdifferenceequ}
\frac {C_3x_n}{\log s_n}= 2\log (n+M), \ \text{for}\ n\ge1.
\end{equation}
We suppress the dependence of this sequence on $M$ in the sequel.
For the sequence $\{x_n\}_{n=0}^\infty$  satisfying \eqref{generatingdifferenceequ}, it follows from \eqref{cond-L}
that
\begin{equation}\label{almostsurereduction}
P_{x_0}^L(L_n^\gamma\le s_n \ \text{for all}\ n)\ge 1-\sum_{n=1}^\infty\frac1{(n+M)^2}.
\end{equation}
Now \eqref{generatingdifferenceequ} is a difference equation corresponding to the differential equation
\begin{equation}\label{diffequ}
\frac{C_3S'(t)}{\log S(t)}= 2\log (t+M),\ t\ge 1.
\end{equation}
Integrating, it follows that
\begin{equation}\label{S(t)est}
\frac{C_3S(t)}{\log S(t)}\le 2(t+M)(\log( t+M)-1)+c,
\end{equation}
for some constant $c$.
If one substitutes $\frac3{C_3}(t+M)\big(\log (t+M)\big)^2$ for $S(t)$ in the left hand side of \eqref{S(t)est},
one finds that the resulting expression is larger than the right hand side of \eqref{S(t)est} for large $t$.
Since the left hand side of \eqref{S(t)est} is increasing as a function of  $S(t)$, (for $S(t)\ge e$),
it follows that $S(t)\le \frac3{C_3}(t+M)\big(\log( t+M)\big)^2$, for sufficiently large $t$.
It then follows that the solution $\{s_n\}_{n=0}^\infty$ to \eqref{generatingdifferenceequ} satisfies
\begin{equation}\label{sn-est}
s_n\le \hat s_n\equiv C(n+M)\big(\log( n+M)\big)^2,\ n\ge 1,
\end{equation}
for some $C>0$.

Now from \eqref{almostsurereduction} and \eqref{sn-est} we can conclude that with probability at least $1-\sum_{n=1}^\infty\frac1{(n+M)^2}$,
we have $L_n^\gamma\le C(n+M)\big(\log(n+M)\big)^2$, for all $n$. However, this is not good enough to prove recurrence when $k=1$.
In fact though,
from \eqref{cond-L}, \eqref{almostsurereduction} and \eqref{sn-est}, we conclude that with probability at least
$1-\sum_{n=1}^\infty\frac1{(n+M)^2}$,
 $\{L^\gamma_n\}_{n=1}^\infty$ is no larger than
$\{\hat L_n\}_{n=1}^\infty$, where $\hat L_n^\gamma=x_0+\sum_{i=1}^n \hat Z_i$,
and $\{\hat Z_n\}_{n=1}^\infty$ is a sequence of independent random variables with
$\hat Z_n\sim\text{Exp}(\frac{C_3}{\log \hat s_n})$.
By Kolmogorov's strong law and \eqref{sn-est}, it follows that
$\frac{\hat L_n}{\sum_{j=1}^n\log \hat s_j}$ almost surely converges to $\frac1{C_3}$.
Using this with \eqref{sn-est}, it follows that $\{\hat  L_n\}_{n=1}^\infty$ almost surely grows on the order
$n\log n$. Consequently, with probability at least  $1-\sum_{n=1}^\infty\frac1{(n+M)^2}$,
 $\{L^\gamma_n\}_{n=0}^\infty$ grows on an order no  larger than $n\log n$.
Using this with \eqref{H-est-driftdisting}, it follows that
$\sum_{n=0}^\infty H(L^\gamma_n)=\infty$, with probability
at least $1-\sum_{n=1}^\infty\frac1{(n+M)^2}$. But as $M$ is arbitrary, we conclude that this occurs
with probability one,  and thus by Proposition \ref{applicable-t-r}, we conclude
that the diffusion with the two-phase drift is recurrent.

We now assume that $k>1$ and prove transience. Chose $x_0$
 sufficiently large so that
the bound in \eqref{expparameter} holds for $x\ge x_0$.
From \eqref{expparameter} and \eqref{transpd}, we have
\begin{equation}\label{cond-Lagain}
P_{x_0}^L(L^\gamma_{n+1}-L_n^\gamma\ge  x|L^\gamma_n)\ge\exp\big(-\frac {C_4x}{(\log L^\gamma_n)(\log^{(2)}L_n^\gamma)^{k-1}}\big).
\end{equation}
Thus, it follows from the law of large numbers that
\begin{equation}\label{1over-e}
\lim_{n\to\infty}\frac1n\{j\le n: L^\gamma_{j+1}-L_j^\gamma\ge \frac {(\log L^\gamma_j)(\log^{(2)}L_j^\gamma)^{k-1}}{C_4}\}\ge e^{-1}, \ \text{a.s.}
\end{equation}

The above inequality states that, asymptotically, at least the fraction $\frac1e$ of the increments  $L^\gamma_{j+1}-L_j^\gamma$ will be of size at least
   $ \frac {(\log L^\gamma_j)(\log^{(2)}L_j^\gamma)^{k-1}}{C_4}$. As such, it provides a lower bound on the growth rate of
   $\{L_n^\gamma\}_{n=0}^\infty$. Since the function $\frac {(\log y)(\log^{(2)}y)^{k-1}}{C_4}$ is increasing in $y$, the ``worst case''
   scenario resulting in the least growth would occur if out of the first $n$ increments, the first $[\frac ne]$ increments satisfied
   the above condition, and the rest did not. We can thus get a lower bound on the growth rate as follows.
Making sure that $x_0>e$, so that $\log^{(2)}x_0>0$,
let $\{x_n\}_{n=0}^\infty$ be the sequence  defined by
\begin{equation}\label{generatingdifferenceequagain}
x_{n+1}=\frac{(\log s_n)(\log^{(2)}s_n)^{k-1}}{C_4}, n\ge0,
\end{equation}
where $s_n=\sum_{j=0}^nx_j$.
 Then it follows from
\eqref{1over-e} that
\begin{equation}\label{lowerbound}
P_{x_0}^\gamma(L^\gamma_n\ge s_{[\frac n{2e}]} \ \text{for all  large}\ n)=1.
\end{equation}

As we did in the recurrent case, we analyze the growth rate of $\{s_n\}_{n=0}^\infty$ by looking at the growth rate of the differential
equation associated with the above difference equation for $\{s_n\}_{n=0}^\infty$.
The differential equation is
$$
\frac{S'(t)}{(\log S(t))(\log^{(2)}S(t))^{k-1}}=\frac1{C_4}.
$$
Integrating this, the leading order  term on the left is $\frac{S(t)}{(\log S(t))(\log^{(2)}S(t))^{k-1}}$, and thus, for large $S(t)$, we have
\begin{equation}\label{S(t)estagain}
\frac{2S(t)}{(\log S(t))(\log^{(2)}S(t))^{k-1}}\ge\frac1 {C_4}t+c,
\end{equation}
for some constant $c$.
If one substitutes $\frac1{2C_4}t(\log t)(\log^{(2)}t)^{k-1}$ for $S(t)$ in \eqref{S(t)estagain}, one finds that the resulting expression
is smaller than the right hand side of \eqref{S(t)estagain} for large $t$. Since the left hand side of \eqref{S(t)estagain}
is increasing as a function of $S(t)$ (for $S(t)$ sufficiently large), it follows that
$S(t)\ge \frac1{2C_4}t(\log t)(\log^{(2)}t)^{k-1}$,  for sufficiently large $t$.
It then follows that the solution $\{s_n\}_{n=0}^\infty$ to \eqref{generatingdifferenceequagain} satisfies
\begin{equation}\label{sn-estagain}
s_n\ge  Cn(\log n)(\log^{(2)}n)^{k-1}, \ n\ge 3,
\end{equation}
for some $C>0$.
We conclude from \eqref{lowerbound} and \eqref{sn-estagain} that
\begin{equation}\label{lowerbound-final}
P_{x_0}^\gamma(L^\gamma_n\ge C_0n(\log n)(\log^{(2)}n)^{k-1} \ \text{for all  large}\ n)=1,
\end{equation}
for some $C_0>0$.
From \eqref{lowerbound-final} and \eqref{H-est-driftdisting}, it follows that
$\sum_{n=0}^\infty H(L^\gamma_n)<\infty$ a.s., and thus by Proposition \ref{applicable-t-r}, we conclude
that the diffusion with the two-phase drift is transient.
\hfill $\square$

\section{Proof of Theorem \ref{down-crossingdisting}}
By comparison, it suffices to consider the case that
$\gamma(x)=\frac1{2b}\log^{(2)}x+\frac k{2b}\log^{(3)}x$,
for $x\ge x_0$, with $x_0$ large enough so that $\log^{(3)}x_0$ is defined.
We need to show   recurrence in the  case that $k=1$, and transience in the case
that $k>1$.

Since $b^T\equiv b$, we have similar to
 \eqref{uTbconst},
 \begin{equation}\label{uTbgamma}
\big(u_T(s)-u_T(s-\gamma(s))\big)\exp(\int_{z_0}^{s-\gamma(s)}2b^T(y)dy)=\frac1{2b}(1-\exp(-2b\gamma(s))).
\end{equation}
Since $b^R=0$, we have from \eqref{uRfunc} that $u^R(x)=x-z_0$.
Thus, from \eqref{applicablecond} we have
\begin{equation}\label{H2}
H(s)=\frac{1-\exp(-2b\gamma(s))}{2b(s-z_0-\gamma(s))+1-\exp(-2b\gamma(s))}.
\end{equation}
Since $\gamma(s)=o(s)$, we conclude from \eqref{H2} that there exist constants $C_1, C_2>0$ such that
\begin{equation}\label{H-est-down-crossingdisting}
\frac{C_1}s\le H(s)\le \frac{C_2}s,\ \text{for large}\ s.
\end{equation}

We now investigate the growth rate of the Markov process $\{L_n^\gamma\}_{n=0}^\infty$. Recall that given $L_j^\gamma=x$, the distribution
of $L_{j+1}^\gamma-L_j^\gamma$
 is the distribution given in \eqref{transpd}.
Since $b^T\equiv b$, we have from \eqref{uTfunc} that $u_T(x)=\frac1{2b}\big(1-\exp(-2b(x-z_0))\big)$.
Thus,
\begin{equation*}
\frac{u'_T(z)}{u_T(z)-u_T(z-\gamma(z))}=\frac{2b}{\exp(2b\gamma(z))-1}.
\end{equation*}
Plugging into this equation the formula for $\gamma(z)$ given above, we have
\begin{equation}\label{expfunc}
\frac{u'_T(z)}{u_T(z)-u_T(z-\gamma(z))}=\frac{2b}{(\log z)(\log^{(2)}z)^k-1}.
\end{equation}
It was shown in the proof of Theorem \ref{driftdisting} that if \eqref{expparameter} holds, then
$\{L_n^\gamma\}_{n=0}^\infty$  grows at least on the  order  $n\log n(\log^{(2)}n)^{k-1}$.
Thus, comparing \eqref{expfunc} with \eqref{expparameter}, it follows that in the case at hand
$\{L_n^\gamma\}_{n=0}^\infty$  grows at least on the order  $n\log n(\log^{(2)}n)^k$.
The same method of proof used to prove that if \eqref{expparameter} holds with $k=1$, then $\{L_n^\gamma\}_{n=0}^\infty$
grows on an order no larger than $n\log n$, also shows that if \eqref{expparameter} holds with $k>1$, then
$\{L_n^\gamma\}_{n=0}^\infty$  grows on an order no larger than $n\log n(\log^{(2)}n)^{k-1}$.
Thus, again comparing \eqref{expfunc} with \eqref{expparameter}, it follows that in the case at hand
$\{L_n^\gamma\}_{n=0}^\infty$  grows on an order no larger than $n\log n(\log^{(2)}n)^k$.
 We conclude that $\{L_n^\gamma\}_{n=0}^\infty$ grows exactly on the order $n\log n(\log^{(2)}n)^k$.
 Using this with \eqref{H-est-down-crossingdisting}, it follow from Proposition \ref{applicable-t-r} that
the diffusion with the two-phase drift is recurrent if $k=1$ and transient if $k>1$.
\hfill $\square$
\section{Proof of Theorem \ref{ballistic}}
Recall that the first $\gamma$-down-crossed time for the process $X(t)$ is given by
$$
\begin{aligned}
&\sigma_\gamma=\inf\{t\ge0: \exists s<t \ \text{with}\ X(t)\le X(s)-\gamma(X(s))\}=\\
&\inf\{t\ge0:  X(t)= X^*(t)-\gamma(X^*(t))\},
\end{aligned}
$$
and  is a stopping time.
Recall that $X^*(\sigma_\gamma)$ has been denoted by $L^\gamma$ and that
$$
\hat\tau_{L^\gamma}=\inf\{t\ge0: X(\sigma_\gamma+t)=L^\gamma\}
$$
is the first time after $\sigma_\gamma$ that the process $X(\cdot)$  returns to its running maximum $L^\gamma$.
Thus, the process $X(\cdot)$ increases from $x_0$ to $L^\gamma$ from time 0 to time
$\sigma_\gamma+\hat\tau_{L^\gamma}$. At the regeneration point $L^\gamma$ at time $\sigma_\gamma+\hat\tau_{L^\gamma}$,
everything begins anew according to the same rules, and also, according to the same distribution,  since the two phases of the drift are constants
and thus independent of location.
If follows from this and the law of large numbers, and the standard technique to go from stopping times to deterministic times,
that
\begin{equation}\label{ballisticequ}
\lim_{t\to\infty}\frac{X(t)}t=\frac{\mathcal{E}_{x_0}L^\gamma-x_0}{\mathcal{E}_{x_0}(\sigma_\gamma+\hat\tau_{L^\gamma})}\ \text{a.s.}
\end{equation}

Recall from the construction in section 3 that  $X(t)$ is in the $Y$-mode up until time $\sigma_\gamma$. Then from time $\sigma_\gamma$ until
time $\sigma_\gamma+\hat\tau_{L^\gamma}$  it is in the $Z$-mode.
Under the assumption of the theorem,
the $Y$-mode corresponds to Brownian motion with a constant drift $b$.
It follows from Doob's optional stopping theorem that
$X(t\wedge\sigma_\gamma)-b(t\wedge\sigma_\gamma)$ is a martingale. Taking expectations, we obtain
\begin{equation}\label{mgexp}
\mathcal{E}_{x_0}X(\sigma_\gamma\wedge t)=x_0+b~\mathcal{E}_{x_0}\sigma_\gamma\wedge t.
\end{equation}
There exists a constant $c>0$ such that a Brownian motion with constant drift $b$ and starting from any $x$ has probability $c$ of downcrossing
the interval $[x-\gamma,x]$ within one unit of time. Thus, it follows that there exists a constant $c_0\in(0,1)$ such
that
\begin{equation}\label{geombd}
\mathcal{P}_{x_0}(\sigma_\gamma>t)\le c_0^t,\ \text{for}\ t\ge 1.
\end{equation}
We have
\begin{equation}\label{c-s}
\mathcal{E}_{x_0}(X(\sigma_\gamma\wedge t);\sigma_\gamma>t)=\mathcal{E}_{x_0}(X(t);\sigma_\gamma>t)
\le(\mathcal{E}_{x_0}X^2(t))^\frac12(\mathcal{P}_{x_0}(\sigma_\gamma>t))^\frac12.
\end{equation}
By comparison, $X(t)$ under $\mathcal{P}_{x_0}$ is stochastically dominated by $x_0+\sqrt aW(t)+bt$, where $W$ is a standard Brownian motion;
thus,
\begin{equation}\label{squareest}
\mathcal{E}_{x_0}X^2(t)\le x_0^2+at+b^2t+2x_0bt.
\end{equation}
Letting $t\to\infty$ in \eqref{mgexp}, and using \eqref{geombd}-\eqref{squareest}, we obtain
\begin{equation}\label{Doob}
\mathcal{E}_{x_0}X(\sigma_\gamma)=x_0+b~\mathcal{E}_{x_0}\sigma_\gamma.
\end{equation}

In the case at hand, where the diffusion coefficient is $a$ instead of 1, the function $u_T$ is given by
$\int_{z_0}^x\exp(-\int_{z_0}^y\frac{2b}adt)=\frac a{2b}\big(1-\exp(-\frac{2b}a(x-z_0))\big)$.
Similar to  \eqref{uTbconst} and \eqref{dbgamma}, it follows that under $\mathcal{E}_{x_0}$,
$L^\gamma-x_0$ is distributed according to an exponential distribution with parameter $\frac{2b}{a\exp(\frac{2 b\gamma}a)-1}$;
thus
\begin{equation}\label{expdist}
\mathcal{E}_{x_0}L^\gamma=x_0+\frac{a\big(\exp(\frac{2 b\gamma}a)-1\big)}{2b}.
\end{equation}
Using \eqref{Doob} and \eqref{expdist}, along with the   fact that $L^\gamma=X(\sigma_\gamma)+\gamma$,
we have
\begin{equation}\label{exp-sigma}
\mathcal{E}_{x_0}\sigma_\gamma=\frac1b\Big(\frac{a\big(\exp(\frac{2 b\gamma}a)-1\big)}{2b}-\gamma\Big).
\end{equation}

We now evaluate $\mathcal{E}_{x_0}\hat\tau_{L^\gamma}$.
From the definition of the process and the fact that the two drift phases are constants and thus independent of location, it follows that
 under $\mathcal{E}_{x_0}$, the distribution of $\hat\tau_{L^\gamma}$ is the distribution
of the first hitting time of $\gamma$ by the diffusion process starting at 0 and corresponding
to the operator $L^{a;c,0,b}\equiv\frac12a\frac{d^2}{dx^2}+b^{c,0,b}(x)\frac d{dx}$,
where $b^{c,0,b}(x)$ is equal to $b$ when $x>0$ and is equal to $c$ when $x\le 0$.
Let $E_0^{a;c,0,b}$ denote the expectation for this diffusion starting from 0.
So
\begin{equation}\label{exp-connection}
\mathcal{E}_{x_0}\hat\tau_{L^\gamma}=E_0^{a;c,0,b}\tau_\gamma.
\end{equation}

We have $E_0^{a;c,0,b}\tau_\gamma=\lim_{N\to\infty}E_0^{a;c,0,b}\tau_\gamma\wedge\tau_{-N}$.
As is well-known \cite{RP95}, \newline $E_0^{a;c,0,b}\tau_\gamma\wedge\tau_{-N}=v_N(0)$,
where $v_N$ solves the equation
$$
\begin{aligned}
&L^{a;c,0,b}v_N=-1 \ \text{in}\ (-N,0)\cup(0,\gamma);\\
&v_N(-N)=v_N(\gamma)=0;\\
&v_N(0^-)=v_N(0^+),\
v_N'(0^-)=v_N'(0^+).
\end{aligned}
$$
Solving this, we obtain
$$
v_N(y)=\begin{cases}\frac {aA_N}{2c}\Big(\exp(\frac{2 cN}a)-\exp(-\frac{2 cy}a)\Big)-\frac{y+N}c,\ -N\le y\le 0;\\
\frac {aD_N}{2b}\Big(\exp(-\frac{2b\gamma}a)-\exp(-\frac{2by}a)\Big)+\frac{\gamma- y}b, \ 0\le y\le \gamma,
\end{cases}
$$
where
$$
\begin{aligned}
&A_N=\frac{\frac\gamma b+\frac Nc+\frac{a(b-c)}{2b^2c}(1-\exp(-\frac{2b\gamma}a))}{\frac a{2c}\big(\exp(\frac{2cN}a)-1\big)+
\frac a{2b}\big(1-\exp(-\frac{2b\gamma}a)\big)};\\
&D_N=A_N+\frac1b-\frac1c.
\end{aligned}
$$
From this we obtain
$$
E_0^{a;c,0,b}\tau_\gamma=\lim_{N\to\infty}v_N(0)=\frac\gamma b+\frac{a(b-c)}{2b^2c}\big(1-\exp(-\frac{2b\gamma}a)\big),
$$
and thus from \eqref{exp-connection},
\begin{equation}\label{exptauhat}
\mathcal{E}_{x_0}\hat\tau_{L^\gamma}=\frac\gamma b+\frac{a(b-c)}{2b^2c}\big(1-\exp(-\frac{2b\gamma}a)\big).
\end{equation}
The theorem now follows from \eqref{ballisticequ}, \eqref{expdist}, \eqref{exp-sigma} and \eqref{exptauhat}.
\hfill $\square$

\end{document}